\documentclass[12pt,a4paper]{amsart}

\usepackage{amsfonts}
\usepackage{amsmath}

\newtheorem{theorem}{Theorem}
\newtheorem{remark}{Remark}
\newtheorem{lemma}{Lemma}

\newtheorem{definition}{Definition}

\def\br{\mathbb R}
\def\bc{\mathbb C}

\def\mand{\quad\mbox{and}\quad}
\def\mbx#1{\makebox[#1cm]{}}
\def\codim{\,\mbox{codim}\,}

\begin{document}
\bibliographystyle{alpha}

\title{A remark on Schubert cells and duality of orbits on flag manifolds}
\author[S. Gindikin]{Simon Gindikin}
\address{\hskip-\parindent
        Simon Gindikin\\
        Department of Mathematics\\
        Hill Center, Rutgers University\\
        110 Frelinghysen Road\\
        Piscataway, NJ 08854-8019, U.S.A.}
\email{gindikin@math.rutgers.edu}

\author[T. Matsuki]{Toshihiko Matsuki}
\address{\hskip-\parindent
        Toshihiko Matsuki\\
        Faculty of Integrated Human Studies\\
        Kyoto University\\
        Kyoto 606-8501, Japan}
\email{matsuki@math.h.kyoto-u.ac.jp}

\date{}

\begin{abstract}
It is known that the closure of an arbitrary $K_\bc$-orbit on a flag manifold is expressed as a product of a closed $K_\bc$-orbit and a Schubert cell (\cite{M2}, \cite{Sp}). We already applied this fact to the duality of orbits on flag manifolds (\cite {GM}). We refine here this result and and give its new applications to the study of domains arising from the duality.
\end{abstract}

\maketitle

\section{Duality of orbits on flag manifolds}

Let $G_\bc$ be a connected complex semisimple Lie group and $G_\br$ a connected real form of $G_\bc$. Let $K_\bc$ be the complexification in $G_\bc$ of a maximal compact subgroup $K$ of $G_\br$. Let $X=G_\bc/P$ be a flag manifold of $G_\bc$ where $P$ is an arbitrary parabolic subgroup of $G_\bc$. Then there exists a natural one-to-one correspondence between the set of $K_\bc$-orbits $S$ and the set of $G_\br$-orbits $S'$ on $X$ given by the condition:
\begin{equation}
S\leftrightarrow S'\Longleftrightarrow S\cap S'\mbox{ is non-empty and compact} \tag{A}
\end{equation}
(\cite{M3}). In the following, we will identify orbits $S$ with $K_\bc$-$P$ double cosets and $S'$ with $G_\br$-$P$ cosets.

We defined in \cite{GM} a subset $C(S)$ of $G_\bc$ by
$$C(S)=\{x\in G_\bc\mid xS\cap S'\mbox{ is non-empty and compact in }X=G_\bc/P\}$$
where $S'$ is the $G_\br$-orbit on $X$ given by (A).

If $S$ is closed, then $S'$ is open (\cite{M1}) and so the condition
$$xS\cap S'\mbox{ is non-empty and compact in }G_\bc/P$$
implies
$$xS\subset S'.$$
Hence the set $C(S)_0$ is the cycle domain (cycle space) for $S'$ (\cite{WW}) where $C(S)_0$ denote the connected component of $C(S)$ containing the identity.

On the other hand, let $S_0$ denote the unique open $K_\bc$-$B$ double coset in $G_\bc$ where $B$ is a Borel subgroup of $G_\bc$ contained in $P$. (We will keep this notation for the whole note.) Then $S'_0$ is the only one closed $G_\br$-$B$ double coset in $G_\bc$ and the condition
$$xS_0\cap S'_0\mbox{ is non-empty and compact in }G_\bc/B$$
implies
$$xS_0\supset S'_0.$$
Let $\{S_j \mid j\in J\}$ be the set of $K_\bc$-$B$ double cosets in $G_\bc$ of codimension one and $T_j=S_j^{cl}$ denote the closure of $S_j$. The sets $T_j$ will play an important role in our constructions.

The compliment of $S_0$ in $G_\bc$ is written as
$$\bigcup_{j\in J} T_j$$
(by Theorem 2 in Section 2). So the set $C(S_0)$ is the compliment of the infinite family of complex hypersurfaces
$$T_j^{-1}g\quad (j\in J,\ g\in S'_0)$$
and hence the connected component $C(S_0)_0$ is Stein.

This domain is sometimes called the ``Iwasawa domain'' since it is a maximal domain where all Iwasawa decompositions can be holomorphically extended from $G_\br$. 

In \cite{GM}, we defined
$$C=\bigcap C(S)$$
where we take the intersection for all $K_\bc$-orbits on $X$ on all flag manifolds $X=G_\bc/P$ of $G_\bc$ and conjectured
$$C=\widetilde{D_0}Z$$
(Conjecture 1.3) where $D_0=\widetilde{D_0}/K_\bc$ is the domain introduced by \cite{AG} (which is sometimes denoted as $\Omega_{AG}$) and $Z$ is the center of $G_\bc$. For connected components, it means
\begin{equation}
C_0=\widetilde{D_0}. \tag{B}
\end{equation}

It is proved in Proposition 8.3 of \cite{GM} that
$$C_0=C(S_0)_0.$$
In other words, $C(S)_0$ is minimal when $S=S_0$. We believe that it is one of central facts of this theory since it gives a very strong estimate of all $C(S)$ through $C(S_0)$ only. So the conjecture (B) is equivalent to
$$C(S_0)_0=\widetilde{D_0}.$$
Recently, the inclusion $C(S_0)_0\subset \widetilde{D_0}$ is proved by Barchini (\cite{B}). On the other hand, the opposite inclusion $C(S_0)_0\supset \widetilde{D_0}$ is proved for all classical cases (\cite{GM}, \cite{KS}) and exceptional Hermitian cases (\cite{GM}). It is announced that it is proved for all spaces in \cite{H}. (We had no chance to see this preprint.)

\begin{remark} \ {\rm In \cite{FH}, the authors deduce the inclusion $C_0 \subset \widetilde{D_0}$ from their result about $C(S)$ for closed $S$ and Proposition 8.1 in \cite{GM}. As we showed above, this inclusion is already the consequence of Proposition 8.3 in \cite{GM} and \cite{B}. So it does not need the results in \cite{FH}.
}\end{remark}

\section{Schubert cells in the category of $K_\bc$-$B$ double cosets}

The principal idea of our considerations in \cite{GM} was that $C(S)_0$ will be essentially independent of neither $S$ nor the flag manifold $X=G_\bc/P$. To justify it, we need to build bridges between $C(S)$ for different $S$ and for it we need to see connections between different $K_\bc$-orbits. It turns out that Schubert cells are very efficient tool for such considerations as in Section 2 and Section 8 in \cite{GM}. They give a possibility to obtain an important information about general $C(S)$ from a consideration of simplest $S$. Here we refine connections between $K_\bc$-orbits and Schubert cells and give more examples of applications.

For a simple root $\alpha$ in the root system with respect to the order defined by $B$, we can define a parabolic subgroup
$$P_\alpha=B\cup Bw_\alpha B$$
of $G_\bc$ such that $\dim_\bc P_\alpha=\dim_\bc B+1$.

\begin{lemma} \ Let $S_1$ be a $K_\bc$-$B$ double coset. Then we have$:$

{\rm (i)} \ If $\dim_\bc S_1P_\alpha=\dim_\bc S_1$, then $S_1^{cl}P_\alpha=S_1^{cl}$.

{\rm (ii)} \ If $\dim_\bc S_1P_\alpha=\dim_\bc S_1+1$, then there exists a $K_\bc$-$B$ double coset $S_2$ such that $S_1^{cl}P_\alpha=S_2^{cl}$.
\end{lemma}

Proof. Though this lemma follows easily from \cite{M2} Lemma 3, we will give a proof for the sake of completeness. Write $S_1=K_\bc gB$. Then we have a natural bijection
$$(g^{-1}K_\bc g\cap P_\alpha)\backslash P_\alpha/B\cong K_\bc\backslash K_\bc gP_\alpha/B=K_\bc\backslash S_1P_\alpha/B$$
by the map $x\mapsto gx$.

(i) \ If $\dim_\bc S_1P_\alpha=\dim_\bc S_1$, then $(g^{-1}K_\bc g\cap P_\alpha)B/B$ is  Zariski open in $P_\alpha/B=P^1(\bc)$ and hence it is dense. So we have
$$S_1^{cl}=(K_\bc gB)^{cl}\supset S_1P_\alpha\supset S_1$$
and therefore $S_1^{cl}=S_1^{cl}P_\alpha$.

(ii) \ Suppose $\dim_\bc S_1P_\alpha=\dim_\bc S_1+1$. Then there exists a $p\in P_\alpha$ such that $(g^{-1}K_\bc g\cap P_\alpha)pB/B$ is Zariski open in $P_\alpha/B=P^1(\bc)$ since the number of $K_\bc$-$B$ double cosets in $G_\bc$ is finite. If we write $S_2=K_\bc gpB$, then we have
$$(S_2)^{cl}\supset S_1P_\alpha\supset S_2$$
and therefore $S_2^{cl}=S_1^{cl}P_\alpha$. \hfill q.e.d.

\begin{theorem} \ Let $S_1$ be a $K_\bc$-$B$ double coset in $G_\bc$ and $w$ an element of $W$. Then we have$:$

{\rm (i)} \ $S_1^{cl}(BwB)^{cl}=S_2^{cl}$ for some $K_\bc$-$B$ double coset $S_2$.

{\rm (ii)} \ $($minimal expression$)$ There exists a $w'\in W$ such that $w'<w$ $($Bruhat order$)$, \ $\ell(w')=\dim_\bc S_2 -\dim_\bc S_1$ and that
$$S_1^{cl}(Bw'B)^{cl}=S_2^{cl}.$$
Here $\ell(w')=\dim_\bc Bw'B-\dim_\bc B$ is the length of $w'$.
\end{theorem}

Proof. (i) \ This follows from Lemma 1 because every Schubert cell $(BwB)^{cl}$ is written as
$$(BwB)^{cl}=P_{\alpha_1}\cdots P_{\alpha_\ell}$$
where $w=w_{\alpha_1}\cdots w_{\alpha_\ell}$ is a minimal expression of $w\in W$.

(ii) \ By Lemma 1, we can choose a subsequence $\beta_1,\ldots,\beta_q$ ($q=\dim_\bc S_2 -\dim_\bc S_1$) of $\alpha_1,\ldots,\alpha_\ell$ such that
$$\dim_\bc S_1^{cl}P_{\beta_1}\cdots P_{\beta_k}=\dim_\bc S_1^{cl}P_{\beta_1}\cdots P_{\beta_{k-1}} + 1$$
for $k=1,\ldots,q$ and that
$$S_2^{cl} =S_1^{cl}(BwB)^{cl} =S_1^{cl}P_{\alpha_1}\cdots P_{\alpha_\ell} =S_1^{cl}P_{\beta_1}\cdots P_{\beta_q} =S_1^{cl}(Bw'B)^{cl}$$
with $w'=w_{\beta_1}\cdots w_{\beta_q}$. \hfill q.e.d.

\begin{remark} {\rm \ $S_1^{cl}(BwB)^{cl}=S_2^{cl}$ implies $S_1^{cl}\subset S_2^{cl}$. But $S_1^{cl}\subset S_2^{cl}$ does not always imply $S_1^{cl}(BwB)^{cl}=S_2^{cl}$ for some $w$ (c.f. \cite{M2}).
}\end{remark}

\begin{definition} \ {\rm For every $K_\bc$-$B$ double coset $S$, we can define, by Theorem 1, a subset $J(S)$ of $J$ by
$$J(S)=\{j\in J\mid S^{cl}(BwB)^{cl}=T_j\mbox{ for some }w\in W\}.$$
}\end{definition}

\begin{lemma} \ Let $S$ be a non-open $K_\bc$-$B$ double coset. Then there exists a simple root $\alpha$ such that
$$\dim_\bc SP_\alpha=\dim_\bc S +1.$$
\end{lemma}

Proof. Write $G_\bc=(Bw_0 B)^{cl}=P_{\alpha_1} \cdots P_{\alpha_m}$ with the longest element $w_0$ in $W$. If
$$\dim_\bc SP_\alpha=\dim_\bc S$$
for all simple roots $\alpha$, then we have, by Lemma 1,
$$G_\bc=S^{cl} G_\bc=S^{cl} P_{\alpha_1} \cdots P_{\alpha_m}=S^{cl},$$
a contradiction. \hfill q.e.d.

\begin{theorem} \ If $\ell(w)<{\rm codim}_\bc S$, then
$$S^{cl}(BwB)^{cl}\subset T_j$$
for some $j\in J(S)$.
\end{theorem}

Proof. Since ${\rm codim}_\bc S^{cl}(BwB)^{cl}=d>0$, we can choose simple roots $\alpha_1,\ldots,\alpha_{d-1}$ such that
$${\rm codim}_\bc S^{cl}(BwB)^{cl}P_{\alpha_1}\cdots P_{\alpha_{d-1}}=1$$
by Lemma 2. Since $(BwB)^{cl}P_{\alpha_1}\cdots P_{\alpha_{d-1}}=(Bw'B)^{cl}$ for some $w'\in W$, we have
$$S^{cl}(BwB)^{cl}\subset S^{cl}(Bw'B)^{cl}=T_j$$
for some $j\in J(S)$. \hfill q.e.d.

\section{Applications}

\begin{definition} \ {\rm For every subset $J'$ in $J$, we define a domain $\Omega(J')$ in $G_\bc$ by
$$\Omega(J')=\{x\in G_\bc\mid xT_j\cap S'_0=\phi\mbox{ for all }j\in J'\}_0.$$
}\end{definition}

We can prove the following corollary:

\bigskip
\noindent {\bf Corollary}\quad {\it Let $S$ be a closed $K_\bc$-$P$ double coset in $G_\bc$. Write $S=S_1^{cl}$ with the dense $K_\bc$-$B$ double coset $S_1$ in $S$. Then we have
$$C(S)_0=\Omega(J(S_1)).$$
}

\begin{remark} \ {\rm (i) \ We can see $C(S_0)_0=\Omega(J)$. By the same argument as for $C(S_0)_0$ in Section 1, we can prove $\Omega(J')$ is Stein for every subset $J'$ in $J$. So the Steinness of $C(S)_0$ (\cite{W}) becomes a corollary of this equivalence $C(S)_0=\Omega(J(S_1))$ (c.f. \cite{HW}).

(ii) \ It is clear that $\Omega(J')\supset \Omega(J)$ for every subset $J'$ in $J$. So we have
$$C(S)_0\supset C(S_0)_0.$$
But this inclusion was already proved in Proposition 8.3 in \cite{GM}. This is natural because the way of proof of the corollary below is essentially the same as that of Proposition 8.3 in \cite{GM}. So the above corollary may be considered as its refinement.
}\end{remark}

Proof of Corollary. Let $x$ be an element on the boundary of $C(S)_0$. Then
$$xS\cap S'_2P\ne\phi$$
for some $G_\br$-$P$ double coset $S'_2P$ in the boundary of $S'$. Here we take $S_2$ as the dense $K_\bc$-$B$ double coset contained in $S_2P$. Since $S$ is right $P$-invariant, we have
$$xS\cap S'_2\ne\phi \mand \dim_\bc S_2>\dim_\bc S.$$

Applying Theorem 1 (ii) to the pair $(S_2^{cl},\ G_\bc)$, we can take a $w\in W$ such that $\ell(w)={\rm codim}_\bc\,S_2$ and that
$$S_2(BwB)^{cl}=G_\bc.$$
So we have $S_2(BwB)^{cl}\supset S_0$ and hence
$$S'_2\subset S'_0(Bw^{-1}B)^{cl}.$$
Since $xS\cap S'_2\ne\phi$, we have
$$xS\cap S'_0(Bw^{-1}B)^{cl}\ne \phi.$$
Hence
$$xS(BwB)^{cl}\cap S'_0\ne\phi$$
which implies $xT_j\cap S'_0\ne\phi$ for some $j\in J(S_1)$ by Theorem 2. Thus $x\notin \Omega(J(S_1))$.

Conversely, suppose
$$xT_j\cap S'_0\ne\phi$$
for some $T_j=S(BwB)^{cl}=S_1^{cl}(BwB)^{cl}$. Note that $j\in J(S_1)$ by Definition 1 and that we may assume $\ell(w)={\rm codim}_\bc \,S -1={\rm codim}_\bc \,S_1 -1$ by Theorem 1 (ii). Then we have
$$xS\cap S'_0(Bw^{-1}B)^{cl}\ne \phi$$
and hence
$$xS\cap S'_3\ne\phi$$
for some $K_\bc$-$B$ double coset $S_3$ such that $S'_3\subset S'_0(Bw^{-1}B)^{cl}$. Hence $S_3(BwB)^{cl}\supset S_0$ and therefore $\dim_\bc S_3\ge \dim_\bc G_\bc -\ell(w)>\dim_\bc S$. So we have
$$S'_3\cap S'=\phi$$
because $S'$ is the union of $G_\br$-$B$ double cosets $S'_4$ satisfying $S_4\subset S$. Hence we have
$$xS\not\subset S'$$
and therefore
$$\mbx6 x\notin C(S).\mbx5 \mbox{q.e.d.}$$

\begin{remark} \ {\rm (i) \ The condition $\ell(w)={\rm codim}_\bc S -1$ does ``not always'' imply
$${\rm codim}_\bc\, S^{cl}(BwB)^{cl}=1.$$
Counter examples exist already for $G_\br=SU(2,1)$.

(ii) \ The construction of the domain $\Omega(J(S_1))$ is essentially equivalent to the construction of ``Schubert domain'' in \cite{HW}. Unfortunately, their basic definition needs a correction and after this correction their proof of Corollary 3.2 corresponding to our Corollary is not complete. We can see that the proof of Corollary using the results in Section 2 is extremely simple. Let us explain the connection between these two constructions introducing notations in \cite{HW}.

\bigskip
Take a Borel subgroup $B_0$ of $G_\bc$ so that $G_\br B_0$ is closed in $G_\bc$. A Borel subgroup $B$ of $G_\bc$ is called an ``Iwasawa Borel subgroup'' if
$$B=g_0B_0g_0^{-1}\quad \mbox{for some }g_0\in G_\br.$$
Let $Z=G_\bc/Q$ be a flag manifold. Then we can take $Q$ so that $Q\supset B_0$. Every Schubert cell $Y$ in $Z$ for $B$ is written as
$$Y=(Bg_0wQ)^{cl}=(g_0B_0wQ)^{cl}$$
with some $w\in W$. Let $C_0$ be a closed $K_\bc$-$Q$ double coset (do not miss with $C_0$ of Section 1!). The ``incidence divisor'' $H_Y$ is written as
$$H_Y =\{g \mid gC_0\cap Y\ne\phi\} =YC_0^{-1}=(g_0B_0wQ)^{cl}C_0^{-1}=(C_0(Qw^{-1}B)^{cl}g_0^{-1})^{-1}.$$

In this point, in \cite{HW}, it is written: ``If ${\rm codim}\,Y\le \dim C_0 +1$, then $H_Y$ is a hypersurface in $G_\bc$.'' It is wrong as we remarked in (i).
 
But if $\codim H_Y=1$, then 
$$H_Y^{-1}=C_0(Qw^{-1}B_0)^{cl}g_0^{-1}=T_jg_0^{-1}$$
for some $j\in J'=J(C_1)$ (where $C_1$ is the dense $K_\bc$-$B_0$ double coset in $C_0$) and $g_0\in G_\br$ by our notation.

So their definition of $y(D)$ should be corrected to
$$y(D)=\{Y=(g_0B_0wQ)^{cl} \mid \codim H_Y=1\}.$$
Here $D$ is the $G_\br$-orbit dual to $C_0$. (The condition $Y\subset Z\setminus D$ follows from $\codim H_Y=1$ because
\begin{align*}
Y\cap D=\phi & \Longleftrightarrow DY^{-1}=D(Qw^{-1}B_0)^{cl}g_0^{-1} \not\ni e \\
& \Longleftrightarrow D(Qw^{-1}B_0)^{cl} \not\ni g_0 \\
& \Longleftrightarrow D(Qw^{-1}B_0)^{cl}\cap G_\br B_0 =\phi \\
& \Longleftrightarrow C_0(Qw^{-1}B_0)^{cl}\cap K_\bc B_0 =\phi \\
& \Longleftrightarrow \codim C_0(Qw^{-1}B_0)^{cl}\ge 1.)
\end{align*}

The Schubert domain is defined as
$$\Omega_S(D)=\left\{G_\bc \setminus \left(\bigcup_{Y\in y(D)} H_Y\right)\right\}_0.$$
This definition is equivalent to ours because
\begin{align*}
g\notin \bigcup_{Y\in y(D)} H_Y & \Longleftrightarrow g^{-1}\notin T_jg_0^{-1}\mbox{ for all }j\in J'\mbox{ and }g_0\in G_\br \\
& \Longleftrightarrow g^{-1}G_\br B_0\cap T_j=\phi \mbox{ for all }j\in J' \\
& \Longleftrightarrow G_\br B_0\cap gT_j=\phi \mbox{ for all }j\in J'.
\end{align*}

(iii) \ So in their Corollary 3.2 in \cite{HW}, (2) $\codim_Z(Y)=q+1$ does not imply (3) $H_Y$ is a hypersurface in $\Omega$. Thus the proof of Corollary 3.2 is incomplete.
}\end{remark}

\begin{remark} {\rm The problem of the description of the domain of cycles $C(S)_0$ for groups $G_\br$ of Hermitian type is simpler than the general case. Firstly, in this case, $D_0=\widetilde{D_0}/K_\bc$ has a very simple description: $D_0\cong G_\br/K\times \overline{G_\br/K}$ (Proposition 2.2 in \cite{GM}). As usual, the equality $C(S)_0=\widetilde{D_0}$ for $S$ ($\leftrightarrow S'$) of nonholomorphic type is reduced to two inclusions. The proof of $C(S)_0 \subset \widetilde{D_0}$ in \cite{WZ1} had a mistake which was corrected in \cite{WZ2}. The opposite inclusion was checked in \cite{WZ1} for classical Hermitian groups. In Proposition 2.4 of \cite{GM}, we gave a very simple proof of this inclusion for arbitrary groups of Hermitian type which is free of case-by-case considerations: the use of Schubert cells makes this fact almost trivial. The note \cite{WZ2} also contains this fact with a proof referred to \cite{HW} (which is incomplete as is explained in Remark 4) but without an appropriate reference on the preceding proof in \cite{GM}. Moreover it asserts a misleading statement that the preprint \cite{GM} does not contain a direct proof.
}\end{remark}

\end{document}